\documentclass[11pt]{article}
\usepackage{fullpage}
\usepackage{amsmath,amsthm,amssymb,amsfonts,graphicx,multicol,subfigure}
\usepackage{cite}
\usepackage[all]{xy}

\def\Z{{\mathbb Z}}

\def\cG{{\mathcal G}}

\newcommand{\dquot}[2]{\raise.5ex\hbox{$#1$}\hbox{\Large$/$}\lower.25ex\hbox{$#2$}}

\newcommand{\bigdquot}[2]{\raise1.5ex\hbox{$#1$}\hbox{\Huge$/$}\lower.75ex\hbox{$#2$}}

\def\Stab{\operatorname{Stab}}
\def\Aut{\operatorname{Aut}}

\usepackage[usenames]{color}

\definecolor{Red}{rgb}{1.0,0.00,0.00}

\theoremstyle{plain}
\newtheorem{thm}{Theorem}[section]

\newtheorem{lma}[thm]{Lemma}

\newtheorem{prop}[thm]{Proposition}

\theoremstyle{definition}
\newtheorem{definition}[thm]{Definition}

\date{
   } 

\begin{document}

\title{Automorphism groups of a graph and a vertex-deleted subgraph\thanks{This work was supported by National Science Foundation grant DMS-0354008.}}

\author{
	Stephen G. Hartke\thanks{This author was supported in part by a Nebraska EPSCoR First Award and by National Science Foundation grant DMS-0914815.}\\ University of Nebraska-Lincoln \\ \texttt{hartke@math.unl.edu}
		\and
	Hannah Kolb\\ Illinois Institute of Technology \\ \texttt{hkolb@iit.edu}
		\and
	Jared Nishikawa\\ Willamette University  \\ \texttt{jnishika@willamette.edu}
		\and
	Derrick Stolee\\University of Nebraska-Lincoln \\ \texttt{s-dstolee1@math.unl.edu}
}

\maketitle

\begin{abstract}
	Understanding the structure of a graph along with the structure of its subgraphs is important for several problems in graph theory.
	Two examples are the Reconstruction Conjecture and isomorph-free generation.	
	This paper raises the question of which pairs of groups can be represented 
		as the automorphism groups of a graph and a vertex-deleted subgraph.  
	This, and more surprisingly the analogous question for edge-deleted subgraphs, 
		are answered in the most positive sense using concrete constructions.
\end{abstract}

\section{Introduction}%

The Reconstruction Conjecture of Ulam and Kelley famously states that the isomorphism class of all graphs graphs on three or more vertices
	is determined by the isomorphism classes of its vertex-deleted subgraphs (see \cite{MR0252258} for a survey of classic results on this problem).
A frequent issue when attacking reconstruction problems is that automorphisms of the substructures lead to ambiguity
	when producing the larger structure.

This paper considers the relation between the automorphism group of a graph and the automorphism groups of
	 the vertex-deleted subgraphs and edge-deleted subgraphs.  
If a group $\Gamma_1$ is the automorphism group of a graph $G$, 
	and another group $\Gamma_2$ is the automorphism group of $G-v$ for some vertex $v$, then
	we say $\Gamma_1$ \emph{deletes to} $\Gamma_2$.
This relation is denoted $\Gamma_1 \to \Gamma_2$.
A corresponding definition for edge deletions is also developed.
Our main result is that any two groups delete to each other, with vertices or edges.

These relations also appear in McKay's isomorph-free generation algorithm \cite{MR1606516}, which is frequently used to enumerate
	 all graph isomorphism classes.
After generating a graph $G$ of order $n$, graphs of order $n+1$ are created by adding vertices and considering each $G+v$.
To prune the search tree, the canonical labeling of $G+v$ is computed, 
	usually by \texttt{nauty}, McKay's canonical labeling algorithm \cite{nauty, HRnauty}.
Finding a canonical labeling of a graph reveals its automorphism group.
Since $G$ was generated by this process, its automorphism group is known
	but is not used while computing the automorphism group of $G+v$.
If some groups could not delete to the automorphism group of $G$, 
	then they certainly cannot appear as the automorphism group of $G+v$ which may allow for some improvement to the canonical labeling algorithm.
The current lack of such optimizations hints that no such restrictions exist, but this notion has not been formalized before this paper.

One reason why this problem has not been answered is that the study of graph symmetry is very restricted, 
	mostly to forms of symmetry requiring vertex transitivity.
These forms of symmetry are useless in the study of the Reconstruction Conjecture, as regular graphs are reconstructible.
On the opposite end of the spectrum, almost all graphs are \emph{rigid} (have trivial automorphism group) \cite{MR1864966}.
Graphs with non-trivial, but non-transitive, automorphisms have received less attention.

Graph reconstruction and automorphism concepts have been presented together before \cite{MR1373683, MR1971819}.
However, there appears to be no results on which pairs of groups allow the deletion relation.
While our result is perhaps unsurprising, it is not trivial.
The reader is challenged to produce an example for $\Z_2 \to \Z_3$ before proceeding.

For notation, $G$ always denotes a graph, while $\Gamma$ refers to a group.  The trivial group $I$ consists of only the
	identity element, $\varepsilon$.  All graphs in this paper are finite, simple, and undirected, unless specified otherwise.
All groups are finite.
The automorphism group of $G$ is denoted $\Aut(G)$ and the stabilizer of a vertex $v$
	in a graph $G$ is denoted $\Stab_{G} (v)$.

\section{Definitions and Basic Tools}
We begin with a formal definition of the deletion relation.

\begin{definition}
	Let $\Gamma_{1},\Gamma_{2}$ be finite groups.  If there exists a graph $G$ with $|V(G)| \geq 3$ and vertex $v \in V(G)$ so that
		$\Aut(G) \cong \Gamma_{1}$ and $\Aut(G-v) \cong \Gamma_{2}$, then $\Gamma_{1}$ (\emph{vertex}) \emph{deletes to} $\Gamma_{2}$,
		denoted $\Gamma_{1}\to \Gamma_{2}$.
	Similarly, the group $\Gamma_1$ \emph{edge deletes to} $\Gamma_2$ if there exists a graph $G$ and edge $e \in E(G)$ so that
		$\Aut(G) \cong \Gamma_1$ and $\Aut(G-e) \cong \Gamma_2$.
	If a specific graph $G$ and subobject $x$ give $\Aut(G) \cong \Gamma_1$ and $\Aut(G-x) \cong \Gamma_2$, the deletion relation
		may be presented as $\Gamma_1 \stackrel{G-x}{\longrightarrow} \Gamma_2$.
\end{definition}

To determine the automorphism structure of a graph, vertices that are not in the same orbit can be distinguished by means of neighboring structures.
A useful gadget to make such distinctions is the rigid tree $T(n)$, where $n$ is an integer at least $2$.  
Build $T(n)$ by starting with a path $u_0,z_1,\dots,z_n$.  
For each $i, 1 \leq i \leq n$, add a path $z_i,x_{i,1},x_{i,2},\dots,x_{i,2i},u_i$ of length $2i+1$.  
This results in a tree with $n+1$ leaves.  
Note that each leaf $u_i$ is distance $2i+1$ to a vertex of degree 3 (except for $u_n$, which is distance $2n+2$).
Thus, the leaves are in disjoint orbits and $T(n)$ is rigid.  
Also, if any leaf $u_i$ is selected with $i \geq 1$, $T(n)-u_i$ is rigid.  
This gives an example of the deletion relation $I \to I$.
For notation, let $J$ be a set and $\{T_j\}_{j\in J}$ be disjoint copies of $T(n)$.  
Then $u_i(T_j)$ designates the copy of $u_i$ in $T_j$.  
This is well-defined since there is a unique isomorphism between each $T_j$ and $T(n)$.

For any group $\Gamma$, a simple, unlabeled, undirected graph $G$ exists with $\Aut(G) \cong \Gamma$ \cite{MR1557026}.
The construction is derived from the well-known Cayley graph\footnote{In most uses of the Cayley graph, a generating set is specified. For simplicity, we use the entire group.}. 
Define $C(\Gamma)$ to be a graph with vertex set $\Gamma$
	and complete directed edge set, where the edge $(\gamma,\beta)$ is labeled with $\gamma ^{-1}\beta$, the element whose 
	right-multiplication on $\gamma$ results in $\beta$.  
The automorphism group of $C(\Gamma)$ is $\Gamma$, and the action on the vertices follows right multiplication by elements in $\Gamma$.
That is, if $\gamma \in \Gamma$, the permutation $\sigma_\gamma$ will take a vertex $\alpha$ to the vertex $\alpha\gamma$.

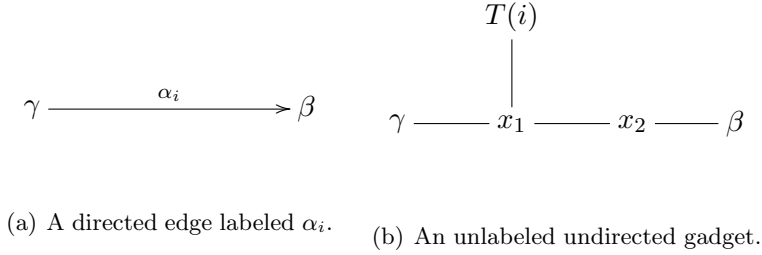
\begin{figure}[htp]
\centering
\mbox
	{
	\subfigure[A directed edge labeled $\alpha_i$.]
		{
			\label{subfig:addflag1}
				$\xymatrix{
				& & & \\
				\gamma \ar[rrr]^{\alpha_i} & \ & \ &  \beta\\
				& & & 
			}$
		}
	\quad	
	\subfigure[An unlabeled undirected gadget.]
		{
			\label{subfig:addflag2}
			$\xymatrix{
				& T(i) \ar@{-}[d] & & \\ 
				\gamma \ar@{-}[r] & x_1  \ar@{-}[r] &x_2  \ar@{-}[r] & \beta\\
				& & & 
			}$
		}
	}

\caption{\label{fig:addflags}Converting a labeled directed edge to an undirected unlabeled gadget.}
\end{figure}

This directed graph with labeled edge sets is converted to an undirected and unlabeled graph
	by swapping the labeled edges with gadgets.  
Specifically, order the elements of $\Gamma = \{ \alpha_1,\dots,\alpha_n\}$ 
	so that $\alpha_1 = \varepsilon$. 
For each edge $(\gamma,\beta)$, subdivide the edge labeled $\alpha_i = \gamma^{-1}\beta$ with vertices $x_1,x_2$, 
and attach a copy $T_{\gamma,\beta}$ of $T(i)$ by identifying $u_0(T_{\gamma,\beta})$ with $x_1$.  Note that $i \geq 2$ in these cases,
since $\alpha_i \neq \varepsilon$. See Figure \ref{fig:addflags} for an example of this process.

Denote this modified graph $C'(\Gamma)$.  We refer to it as \emph{the} Cayley graph of $\Gamma$.
Note that the automorphisms of $C'(\Gamma)$ are uniquely determined by the permutation of the group elements and
preserve the original edge labels, since the trees $T(i)$ identify the label $\alpha_i$ and have a unique isomorphism between copies.
Hence, $\Aut(C'(\Gamma)) \cong \Aut(C(\Gamma)) \cong \Gamma$.

\begin{lma}\label{lma:cayleystable}
	Let $\Gamma$ be a group and $G = C'(\Gamma)$.  
	Then the stabilizer of the identity element $\varepsilon$ (as a vertex in $G$) is trivial.
	That is, $\Stab_G(\varepsilon) \cong I$.
\end{lma}

\begin{proof}
	Every automorphism of $G$ is represented by right-multiplication of $\Gamma$.
	Hence, every automorphism except the identity map will displace $\varepsilon$.
\end{proof}

\section{Deletion Relations with the Trivial Group}

Now that sufficient tools are available, we prove some basic properties. 

\begin{prop}(The Reflexive Property)\label{prop:reflex} For any group $\Gamma,$ $\Gamma \to \Gamma$.\end{prop}

\begin{proof}
	Let $\Gamma$ be non-trivial, as the trivial case has been handled by the rigid tree $T(n)$.
	Let $G$ be the Cayley graph $C'(\Gamma)$.  
	Create a supergraph $G'$ by adding a dominating vertex $v$ with a pendant vertex $u$.
	Now, $u$ is the only vertex of degree $1$, and $v$ is the only vertex adjacent to $u$.  Hence, 
		these two vertices are distinguished in $G'$ from the vertices of $G$.
	Removing $v$ leaves $G$ and the isolated vertex $u$.
	Thus, $\Gamma$ is the automorphism group for both $G'$ and $G'-v$.
\end{proof}

A key part of our final proof relies on the trivial group deleting to any group.

\begin{lma}\label{lma:trivialtoall}
	For all groups $\Gamma$, $I \to \Gamma$.
\end{lma}

\begin{proof}
	Let $G = C'(\Gamma)$.
	Let $n = |\Gamma|$.  Order the group elements of $\Gamma$ as $\alpha_1,\dots, \alpha_n$.  
	Create a supergraph, $G'$, by adding vertices as follows: For each $\alpha_i$, 
		create a copy $T_{\alpha_i}$ of $T(2n)$ and identify $u_0(T_{\alpha_i})$ with the vertex $\alpha_i$ in $G$
		(Here, $2n$ is used to distinguish these copies from the edge gadgets).
	Add a vertex $v$ that is adjacent to
	$u_i(T_{\alpha_i})$ for each $i$.  
	For each $\alpha_i$, the leaf of $T_{\alpha_i}$ adjacent to $v$ distinguishes $\alpha_i$.
	Hence, no automorphisms exist in $G'$.  
	However, $G'-v$ restores all automorphisms $\pi$ from $\Aut(G)$ by mapping $T_{\alpha_i}$ 
	to $T_{\pi(\alpha_i)}$ through the unique isomorphism.
\end{proof}

Note that this proof uses a very special vertex that enforces all vertices to be distinguished.
Before producing examples where deleting a vertex removes symmetry, it may be useful to remark that such a
distinguished vertex cannot be used.

\begin{lma}\label{lma:stabilizer}
	Let $G$ be a graph and $v \in V(G)$.
	Then, automorphisms in $G$ that stabilize $v$ form a subgroup in the automorphism group of $G-v$.
	That is, $\Stab_G(v) \leq \Aut(G-v)$.
\end{lma}

\begin{proof}
	Let $\pi \in \Stab_G(v)$.  The restriction map $\pi|_{G-v}$ is an automorphism of $G-v$.
\end{proof}

The implication of this lemma is removing a vertex with a trivial orbit cannot remove automorphisms.  
However, we can remove all symmetry in a graph using a single vertex deletion.

\begin{lma}\label{lma:grouptrivial}
	For any group $\Gamma$, $\Gamma \to I$.
\end{lma}

\begin{proof}
	Assume $\Gamma \not\cong I$, since the reflexive property handles this case.
	Let $G = C'(\Gamma)$ and $n = |\Gamma|$.
	
	Let $G_1, G_2$ be copies of $G$ with isomorphisms $f_1\!: G \to G_1$ and $f_2\!: G \to G_2$.
	Create a graph $G'$ from these two copies as follows.
	For all elements $\gamma$ in $\Gamma$, create a copy $T_\gamma$ of $T(n)$ 
		and identify $u_0(T_\gamma)$ with $f_1(\gamma)$ and $u_n(T_\gamma)$ with $f_2(\gamma)$.
	Note that $\Aut(G') \cong \Gamma$, since no vertices from $G_1$ can map to $G_2$ 
		from the asymmetry of the $T_\gamma$ subgraphs,
		and any automorphism of $G_1$ extends to exactly one automorphism of $G_2$.
	
	Any automorphism $\pi$ of $G' - f_1(\varepsilon)$ must induce an automorphism $\pi|_{G_2}$ of $G_2$.  
	But the vertices of $G_1$ must then permute similarly (by the definition 
		$\pi(f_1(x)) = f_1 f_2^{-1} \pi f_2(x)$).
	Since $f_1(\varepsilon)$ is not in the image of $\pi$, $\pi$ stabilizes $f_2(\varepsilon)$.
	Lemma \ref{lma:cayleystable} implies $\pi$ must be the identity map.  
	Hence, $\Aut(G'-f_1(\varepsilon)) \cong I$.
\end{proof}

\section{Deletion Relations Between Any Two Groups}
We are sufficiently prepared to construct a graph to reveal the deletion relation for all pairs of groups.

\begin{thm}\label{thm:vertexall}
	If $\Gamma_1$ and $\Gamma_2$ are groups, then $\Gamma_1 \to \Gamma_2$.
\end{thm}

\begin{proof}
	Assume both groups are non-trivial, since Lemmas \ref{lma:trivialtoall} and \ref{lma:grouptrivial} cover these cases.
	Let $G_1 = C'(\Gamma_1)$.
	Then identify $v_1 \in V(G_1)$ as the vertex corresponding to $\varepsilon \in \Gamma_1$.  
	Note that $\Stab_{G_1}(v_1) \cong I$ as in Lemma \ref{lma:cayleystable}.  
	Also by Lemma \ref{lma:trivialtoall}, there exists a graph $G_2$ and vertex $v_2$ so that
		$I \stackrel{G_2-v_2}{\longrightarrow} \Gamma_2$.
	Define $n_i = |\Gamma_i|$.  Order the elements of $\Gamma_1$ as $\alpha_{1,1},\alpha_{1,2},\dots,\alpha_{1,n_1}$
	so that $\alpha_{1,1} = \varepsilon = v_1$.
	
	We collect the necessary properties of $G_1, G_2, v_1, v_2$ before continuing.  
	First, $G_1$ has automorphisms $\Aut(G_1) \cong \Gamma_1$ and $v_1$ is trivially stabilized ($\Stab_{G_1}(v_1) \cong I$).
	Second, $G_2$ is rigid ($\Aut(G_2) \cong I$) but $G_2-v_2$ has automorphisms $\Aut(G_2-v_2) \cong \Gamma_2$.
	The following construction only depends on these requirements.

	Let $H_1,\dots,H_{n_1}$ be copies of $G_2$.
	Construct a graph $G$ by taking the disjoint union of $G_1$, $H_1$, $\dots,$ $H_{n_1}$, and 
		adding edges between $\alpha_{1,i}$ and every vertex of $H_{i}$, for $i = 1,\dots, n_1$.
	Since $\Aut(H_i) \cong I$, the automorphism group of $G$ cannot permute
	the vertices within each $H_i$.  However, the vertices of $G_1$ can permute freely within $\Aut(G_1) \cong \Gamma_1$, 
	since $H_i \cong H_j$ for all $i, j$.  Hence, $\Aut(G) \cong \Gamma_1$.
	
	When the copy of $v_2$ in $H_1$ is deleted from $G$, the automorphisms of $H_1-v_2$ are $\Gamma_2$.  
	However, the vertex $v_1$ of $G_1$ is now distinguished since it is adjacent to a copy of $G_2-v_2$, unlike the other
	elements of $\Gamma_1$ in $G_1$ which are adjacent to a copy of $G_2$.  
	This means the permutations of $G_1$ must stabilize $v_1$.
	Since $\Stab_{G_1}(v_1) = I$, the only permutation allowed on $G_1$ is the identity.  
	However, $H_1-v_2$ has automorphism group $\Gamma_2$.  	Hence, $\Aut(G-v_2) \cong \Gamma_2$.
\end{proof}

\begin{figure}[htp]
\centering
\mbox
	{
	\subfigure[$G$ with $\Aut G \cong \Gamma_1$.]
		{
		  \includegraphics[height=1.5in]{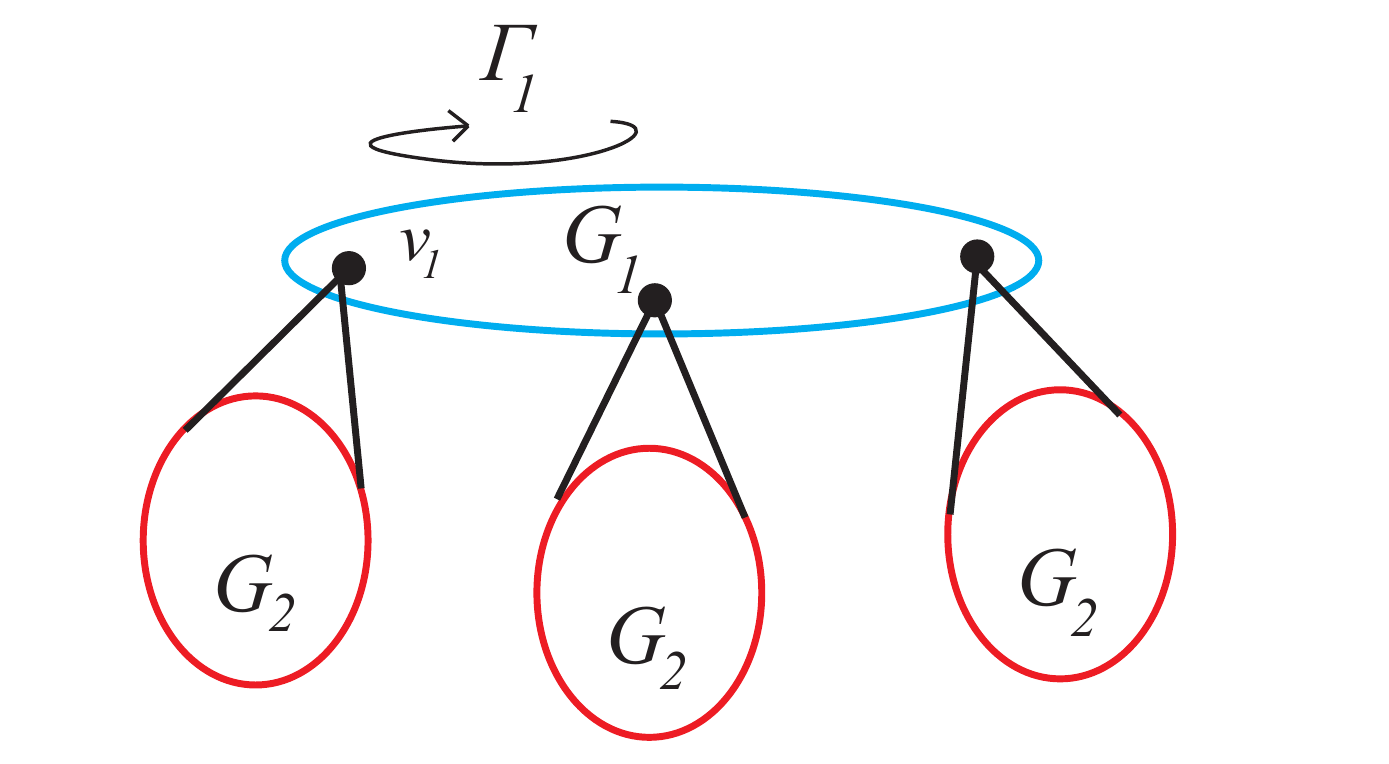}
		}
	\quad	
	\subfigure[$G-v_2$ with $\Aut G-v_2 \cong \Gamma_2$.]
		{
		  \includegraphics[height=1.5in]{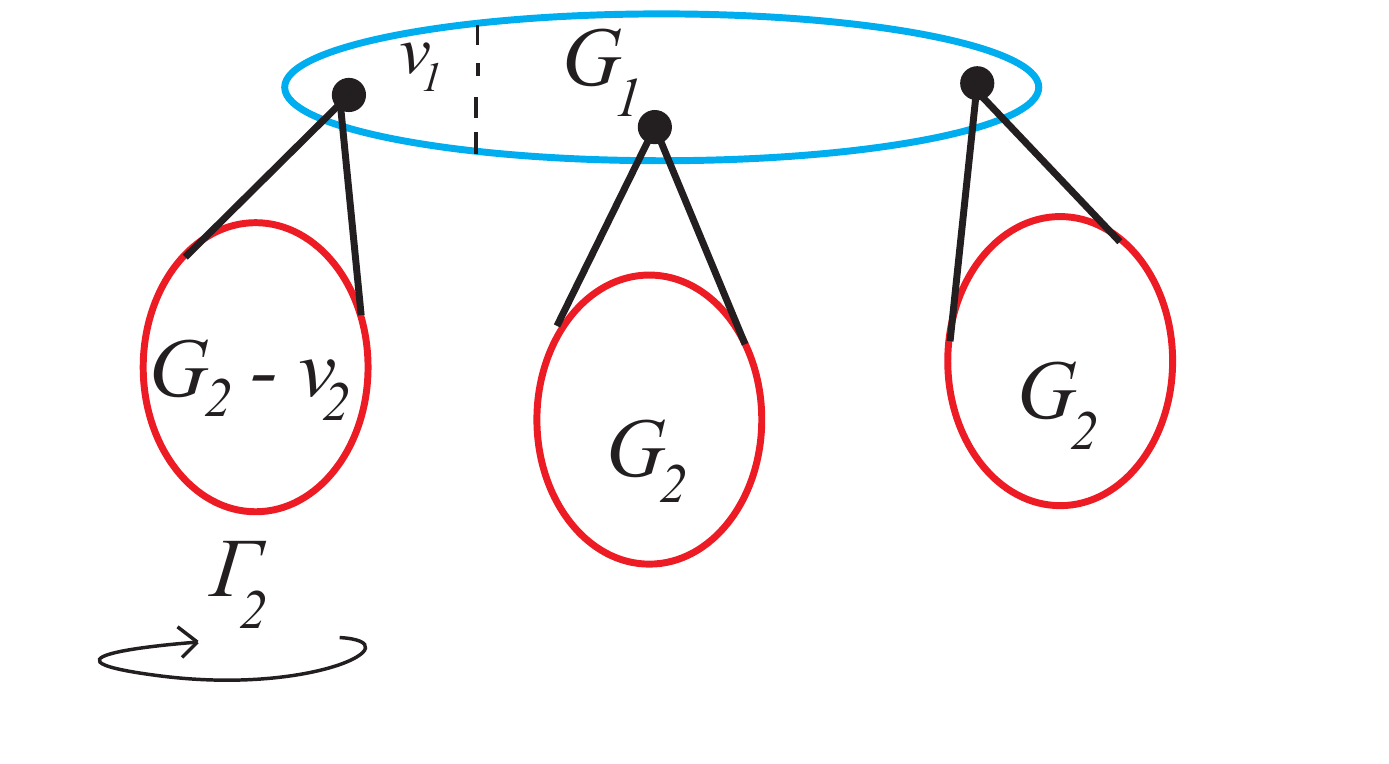}
		}
	}

\caption{\label{fig:vertexall}The vertex deletion construction.}
\end{figure}

Figure \ref{fig:vertexall} presents a visualization of the automorphisms in this construction before and after the deletion.
A very similar construction produces this general result for the edge case.

\begin{thm}\label{thm:edgeall}
	If $\Gamma_1$ and $\Gamma_2$ are groups, then there exists a graph $G$ and an edge $e \in E(G)$ so that
		$\Gamma_1 \stackrel{G-e}{\longrightarrow} \Gamma_2$.
\end{thm}

\begin{proof}
	Set $n_i = |\Gamma_i|$.
	Let $G_1 = C'(\Gamma_1)$ with $v_1$ corresponding to $\varepsilon \in \Gamma_1$ and order the elements of $\Gamma_1$ similarly to the
	proof of Theorem \ref{thm:vertexall}.  
	
	Form $G_2$ by starting with $C'(\Gamma_2)$ and making a copy $T_{\gamma}$ of $T(2n_2)$ for each element $\gamma \in \Gamma_2$, 
	identifying $\gamma \in V(C'(\Gamma_2))$ with $u_0(T_{\gamma})$.  Now, add an edge $e$ between $u_{2n_2}(T_{1})$ and $u_{2n_2-1}(T_1)$.
	This distinguishes the element $\varepsilon$ as a vertex in $C'(\Gamma_2)$ and hence is stabilized.  So, $\Aut(G_2) \cong I$ and if $e$ is removed all
	group elements are symmetric again, so $\Aut(G_2-e) \cong \Gamma_2$.
	
	Notice that $G_1, G_2, v_1, e$ satisfy the requirements of the construction of $G$ in Theorem \ref{thm:vertexall}.  
	Hence, the same construction (with $e$ in place of $v_2$) provides an example of edge deletion from $\Gamma_1$ to $\Gamma_2$.
\end{proof}

Note that the graph produced for Theorem \ref{thm:edgeall} can be used for the proof of Theorem \ref{thm:vertexall} by subdividing $e$ and
using the resulting vertex as the deletion point.


\section{Future Work}
While the question posed in this paper is answered completely for the class of all graphs, there remain questions for special cases.
For instance, the automorphism groups of trees are fully understood \cite{MR607504}. 
Let $\cG_T$ be the class of groups that are represented by the automorphism groups of trees and $\cG_F$ represented by automorphisms of forests\footnote{An elementary proof shows that $\cG_T = \cG_F$.}.
The constructions in this paper are not trees, so new methods will be required to answer the following questions.
If we restrict to trees, can any group in $\cG_T$ delete to any group in $\cG_F$?
Or, if we restrict to deleting leaves (and hence stay connected) can all pairs of groups in $\cG_T$ delete to each other?

\begin{figure}
\centering
%
		  \includegraphics[height=1.5in]{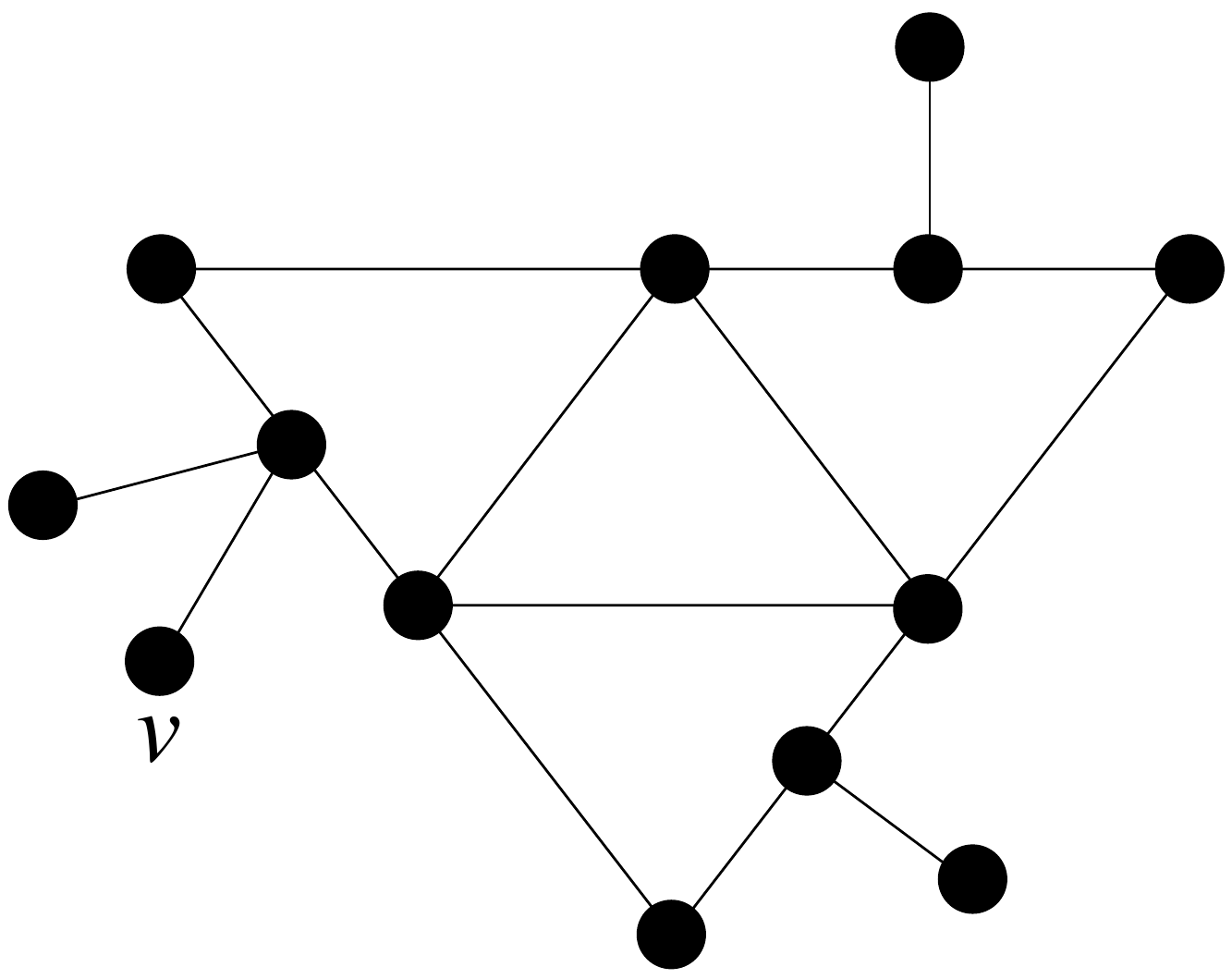}
		  \caption{\label{fig:z2toz3}This graph $G$ has $\Aut(G)\cong \Z_2$ and $\Aut(G-v) \cong \Z_3$.}
\end{figure}

Another interesting aspect of our construction is that the resulting graphs are very large, with the order of the graphs cubic in the size of the groups.
Which of these relations can be realized by small graphs?  Can we restrict the groups that can appear based on the order of the graph?
The current-best upper bound on the order of a graph $G$ with automorphism groups isomorphic to a given group $\Gamma$ is
	$|V(G)| \leq 2|\Gamma|$ and $\Aut(G) \cong \Gamma$ \cite{BabaiOrder}.
This has particular application to McKay's generation algorithm, where only ``small" examples are usually computed 
	(for example, all connected graphs up to 11 vertices were computed in \cite{MR1448235}).  To demonstrate that this
is not trivial, see Figure \ref{fig:z2toz3} for a graph showing $\Z_2\to\Z_3$.

While Theorem \ref{thm:vertexall} shows that there exists a graph where \emph{some} vertex can be deleted to demonstrate the deletion relations, 
	our constructions have many other vertices that behave in very different ways when they are deleted.
When relating to the Reconstruction Conjecture, 
	this raises questions regarding the combinations of automorphism groups that appear in the vertex-deleted subgraphs.
For instance, if the multiset of
vertex-deleted automorphism groups is provided, can one reconstruct the automorphism group?  This question only gives the groups, 
but not the vertex-deleted subgraphs.  An example is that $n$ copies of $S_{n-1}$ must reconstruct to $S_n$, but it is unknown 
whether the graph is $K_n$ or $nK_1$.  Since $\Aut(G) = \Aut(\overline{G})$, this ambiguity will always naturally arise.  
Can it arise in other contexts?  Is the automorphism group recognizable from a vertex deck?

\bibliographystyle{alpha}
\bibliography{bibliography}


\end{document}